\documentclass[11pt]{article}

\usepackage{t1enc,a4,array}
\usepackage[english]{babel}

\usepackage{amsfonts}
\usepackage[dvips]{graphicx}
\usepackage[dvips]{epsfig}
\usepackage{enumerate}
\usepackage{amsthm}
\usepackage{color}

\usepackage{amsmath}
\usepackage{amssymb}
\usepackage{psfrag}
\input xy 
\xyoption{all} 

\newtheorem{theorem}{Theorem}[section]
\newtheorem{definition}[theorem]{Definition}
\newtheorem{lemma}[theorem]{Lemma}

\title{\bf GROUPS ACTING SIMPLY TRANSITIVELY ON HYPERBOLIC TRIANGULAR BUILDINGS}

\author{Lisa Carbone, Riikka Kangaslampi, and Alina Vdovina}

\begin{document}




\maketitle

\begin{abstract}
  We construct and classify all groups, given by triangular
  presentations associated to the smallest thick generalized quadrangle, that
  act simply transitively on the vertices of hyperbolic triangular
  buildings of the smallest non-trivial thickness.  Our classification
  shows 23 non-isomorphic torsion free groups (obtained in an earlier
  work) and 168 non-isomorphic torsion groups acting on one of two
  possible buildings with the smallest thick generalized quadrangle as the link of each vertex.  In analogy with the Euclidean case, we find both torsion and
  torsion free groups acting on the same building.
\end{abstract}


\section{Introduction}

Intensive study of groups acting simply transitively on Euclidean
buildings was initiated in \cite{Cartwright} and
\cite{Cartwright2}. This work has had a considerable impact in several
directions.  For example, this led to new examples of fake projective
planes (\cite{Kato}), and, finally, to their full classification
(\cite{Cartwright3}).  In the Euclidean case, there are two
non-isomorphic buildings of minimal non-trivial thickness admitting a
simply transitive action, and eight isomorphism classes of groups
acting on these two buildings simply transitively and in a type
preserving manner. Within this isomorphism class, five groups are
torsion free and three have torsion \cite{Cartwright2}.

In this paper we study groups acting simply transitively
on hyperbolic buildings with the smallest thick generalized quadrangle as the
link of each vertex.  The torsion free groups acting simply transitively on such buildings
were classified in \cite{KV}.  Here we classify triangle presentations
associated to the smallest thick generalized quadrangle, as well as groups
with torsion coming from these presentations.  

It is known (\cite{Jacek}) that up to isomorphism, there are only two
possible triangular hyperbolic buildings, with the smallest
generalized quadrangle as the link of each vertex, admitting a simply transitive
action. We note that in the formulation of the main theorem in
\cite{Jacek} the appropriate polygonal complexes are required to be
symmetric, but the proof works also for buildings admitting simply
transitive actions.

In \cite{KV} the authors constructed, for any $n$, torsion free groups
acting cocompactly on hyperbolic buildings with $n$-gonal
chambers. Our strategy in this paper is to modify the construction in
\cite{KV} to include the torsion case as well.

Our classification shows 168 non-isomorphic torsion groups acting on
one of two possible buildings with the smallest thick generalized quadrangle
as the link of each vertex.  In analogy with the Euclidean case, we find both torsion
and torsion free groups acting on the same building. These groups are
listed in Appendix 1. The two possible buildings are denoted by (1)
and (2) in Appendix 1.

The link of order 2 for a Kac-Moody building with the minimal
generalized quadrangle as the link of each vertex and equilateral triangular chambers
was computed in unpublished paper by the first author and
D. Cartwright and T. Steger (\cite{Carbone-Cartwright-Steger}), using
an invariant for links of order 2 developed by T. Steger. The
Kac-Moody building coincides with our building with number (2).

By \cite{Jacek} there are only two possible isomorphism classes of
buildings with the smallest thick generalized quadrangle as the link of each vertex and by
results of the present paper at least two of them are
non-isomorphic. Thus all the groups from Appendix 1 with building
number (2) are cocompact  lattices in the automorphism group of the corresponding
Kac-Moody building. It remains to determine if it is possible to embed
these lattices into the corresponding Kac-Moody group.

The existence of cocompact lattices in certain Kac-Moody groups has
already been established. In \cite{CarboneGarland}, the authors
generalized Lubotzky's construction of Schottky groups of
automorphisms in $SL_2$ over a nonarchimedean local field to give
torsion free cocompact lattices in any rank 2 locally compact
Kac-Moody group over a finite field $\mathbb{F}_q$. In \cite{CT}
Capdebosq and Thomas classified cocompact lattices with torsion and
with quotient a simplex in rank 2 Kac-Moody groups corresponding to
symmetric generalized Cartan matrices. In \cite{CarboneCobbs}, the
first author and Cobbs showed that over the field with 2 elements,
rank 3 Kac-Moody groups of noncompact hyperbolic type whose Weyl
groups are a free product of $\mathbb{Z}/2\mathbb{Z}$'s contain a
cocompact lattice that also acts discretely and cocompactly on a
simplicial tree. In \cite{Bourdon} and \cite{Bourdon2}, Bourdon constructed a family of
cocompact lattices in the automorphism groups of certain hyperbolic
Kac-Moody buildings. In \cite{RemyRonan}, R\'emy and
Ronan showed that Bourdon's cocompact lattices $\Gamma_{r,q+1}$,
$r\geq 5,\ q\geq 3$, can be embedded into the closure of right-angled
Kac-Moody groups in the automorphism groups of their buildings,
$I_{r,q+1}$ for $q$ a prime power. 

In all of the above cases, the Kac-Moody buildings are
right-angled. The groups we construct here are the first examples of
cocompact lattices acting on buildings that are not right-angled.

By \cite{Wise} it is known that groups acting cocompactly on
hyperbolic buildings, in such a way that the chamber is a polygon with
at least four sides, are residually finite.  But whether or not groups
acting cocompactly on triangular hyperbolic buildings are residually
finite remains an open question.  Our hyperbolic groups acting simply
transitively on triangular hyperbolic buildings are possible
candidates of such groups that are not residually finite. The commutator
subgroups of many of our examples are perfect groups (that is, they
have trivial abelianizations) and an extensive computer search (which
was carried out since the paper \cite{KV} was completed) did not find
any normal subgroups of these commutator subgroups.

To prove our main theorem, we used a program written in Fortran to determine the equivalence classes of triangular presentations. We used Magma to determine isomorphism classes of dual graphs of polyhedra and hence of triangle presentations.

\section{Definitions and main results}

Recall that a {\em generalized $m$-gon} is a connected, bipartite
graph of diameter $m$ and girth (the length of shortest circuit)
$2m$, in which each vertex lies on at least two edges.

We will call a {\em polyhedron} a two-dimensional complex which is
obtained from several oriented $p$-gons (Euclidean or hyperbolic) with
words on the boundary, by identification of sides with the same labels
respecting orientation. We assume that each side of our polygons has
length 1.

Consider a sphere of a radius $0<\epsilon <1$ at a vertex of the
polyhedron.  The intersection of the sphere with the polyhedron is a
graph, which is called the {\em link} at this point. Consider now
 a sphere of a radius $1+\epsilon$, $0<\epsilon <1$ at a
vertex of the polyhedron. The intersection of this sphere and the
polyhedron will be called a {\em link of order two}.

We will use the definition of a hyperbolic building given in
\cite{Haglund}, where an infinite series of examples of hyperbolic
buildings, with prescribed local structure, were
constructed and studied.

\begin{definition} Let $\mathcal{P}(p,m)$ be a tessellation of the
  hyperbolic plane by regular polygons with $p$ sides, with angles
  $\pi/m$ in each vertex where $m$ is an integer.  A {\em hyperbolic
    building} is a polygonal complex $X$, which can be expressed as
  the union of subcomplexes called apartments such that:

\begin{itemize}
\item[1.] Every apartment is isomorphic to $\mathcal{P}(p,m)$.
\item[2.] For any two polygons of $X$, there is an apartment
containing both of them.
\item[3.] For any two apartments $A_1, A_2 \in X$ containing the same
  polygon, there exists an isomorphism $ A_1 \to A_2$ fixing $A_1 \cap
  A_2$.
\end{itemize}
\end{definition}

Let $C_p$ be a polyhedron whose faces are $p$-gons and whose links are
generalized $m$-gons with $mp>2m+p$. We equip every face of $C_p$ with
the hyperbolic metric such that all sides of the polygons are
geodesics and all angles are $\pi/m$.  Then the universal covering of
such a polyhedron is a hyperbolic building (see \cite{Gaboriau-Paulin}).

Therefore to construct hyperbolic buildings with cocompact group actions, it is
sufficient to construct finite polyhedra with appropriate links.

We recall also the definition of a polygonal presentation introduced
in \cite{Vdovina}:

\begin{definition} Suppose we have $n$ disjoint connected bipartite
  graphs \linebreak $G_1, G_2, \ldots, G_n$.  Let $P_i$ and $L_i$ be the
  sets of black and white vertices respectively in $G_i$,
  $i=1,\dots,n$; let $P=\bigcup P_i$, $L=\bigcup L_i$, $P_i \cap P_j =
  \emptyset$, $L_i \cap L_j = \emptyset$ for $i \neq j$ and let
  $\lambda$ be a bijection $\lambda: P\to Q$.

  A set $\mathcal{K}$ of $k$-tuples $(x_1,x_2, \ldots, x_k)$, $x_i \in
  P$, will be called a {\em polygonal presentation} over $P$
  compatible with $\lambda$ if

  \begin{itemize}
  \item[(1)] $(x_1,x_2,x_3, \ldots ,x_k) \in \mathcal{K}$ implies that
    $(x_2,x_3,\ldots,x_k,x_1) \in \mathcal{K}$;
  \item[(2)] given $x_1,x_2 \in P$, then $(x_1,x_2,x_3, \ldots,x_k)
    \in \mathcal{K}$ for some $x_3,\ldots,x_k$ if and only if $x_2$
    and $\lambda(x_1)$ are incident in some $G_i$;
  \item[(3)] given $x_1,x_2 \in P$, then $(x_1,x_2,x_3, \ldots ,x_k)
    \in \mathcal{K}$ for at most one $x_3 \in P$.
  \end{itemize}

If there exists such $\mathcal{K}$, we will call $\lambda$ a {\em
  basic bijection}.
\end{definition}

{\bf Remark 1.}  The polygonal presentations with $k=3$, $n=1$, and
$G_1$ is the smallest generalized 3-gon have been listed in
\cite{Cartwright} and \cite{Edjvet-Howie}.

We use the following definition of equivalence, which is similar to
the one in \cite{Cartwright2}.

\begin{definition}
  Let $\mathcal{K}_1$ and $\mathcal{K}_2$ be two polygonal
  presentations with $k=3 $, $n=1$, and for which the graph $G_1$ is a
  generalized 4-gon. Then $\mathcal{K }_1$ and $\mathcal{K}_2$ are
  equivalent if there exists an automorphism of the generalized 4-gon
  which transforms the 4-gon of $\mathcal{K}_1$ to the 4-gon of
  $\mathcal{K}_2$.
\end{definition}

Here  we classify all polygonal presentations for $k=3$, $n=1$ and $G_1$ is the smallest thick generalized quadrangle (4-gon). Figure 1 shows the graph $G_1$.

 \begin{figure}
\begin{center}
\includegraphics[scale=0.32]{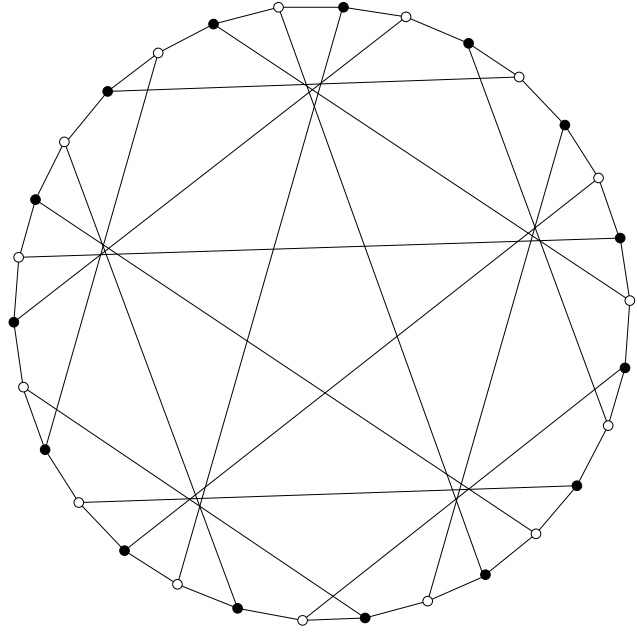}\label{diagram}
\end{center}
\caption{The graph $G_1$}
\end{figure}

In \cite{KV} the authors classified all polygonal presentations for the case
$k=3$, $n=1$ and $G_1$  is the smallest thick generalized quadrangle, when at least two labels in each triangle are different. This
corresponds to the case of torsion free groups acting simply
transitively on the building.

\begin{theorem}(\cite{KV})\label{oldtheorem}
  There are 45 non-equivalent torsion free triangle presentations
  associated to the smallest thick generalized quadrangle. These give
  rise to 23 non-isomorphic torsion free groups, acting simply
  transitively on triangular hyperbolic buildings of smallest
  non-trivial thickness.
\end{theorem}

It turns out that if we allow torsion in the groups acting simply
transitively on hyperbolic triangular buildings, the number of
non-equivalent presentations and the number of non-isomorphic groups
is much larger:

\begin{theorem}\label{maintheorem}
  There are 7159 non-equivalent triangle presentations corresponding
  to groups with torsion associated to the smallest generalized
  quadrangle. These give rise to 168 non-isomorphic groups, acting on
  one of two possible triangular hyperbolic buildings with the
  smallest thick generalized quadrangle as the link of each vertex (listed in Appendix 1).
\end{theorem}

We can associate a polyhedron $X$ on $n$ vertices with each polygonal
presentation $\mathcal{K}$ as follows: for every cyclic $k$-tuple
$(x_1,x_2,x_3,\ldots,x_k)$ we take an oriented $k$-gon on the boundary
of which the word $x_1 x_2 x_3\ldots x_k$ is written. To obtain the
polyhedron we identify the corresponding sides of the polygons,
respecting orientation. We say that the polyhedron $X$ corresponds to
the polygonal presentation $\mathcal{K}$.

The following lemma was proved in \cite{Vdovina}:

\begin{lemma}A polyhedron $X$ which corresponds to a polygonal
  presentation $\mathcal{K}$ has graphs $G_1, G_2, \ldots, G_n$ as
  vertex-links.
\end{lemma}

Polyhedra corresponding to polygonal presentations from Theorem
\ref{oldtheorem} have generalized 4-gons as vertex-links and regular
hyperbolic triangles with angles $\pi/4$ as faces.  The universal
covering of such a polyhedron is a hyperbolic building (see
\cite{Gaboriau-Paulin}).  Moreover, with the metric introduced in
\cite[p.~165]{Ballmann-Brin} this building is a complete metric space of
non-positive curvature in the sense of Alexandrov and Busemann
\cite{Ghys-Harpe}.  Examples of hyperbolic buildings with right-angled
triangles were constructed by M.~Bourdon in \cite{Bourdon} and in
\cite{Gaboriau-Paulin}.

{\bf Remark.} If we have a group with torsion, to apply 
\cite{Gaboriau-Paulin}, we have to consider index 3 subgroups
(obtained in a canonical way by changing alphabets), to form polyhedra with
three vertices, then to go to the universal cover carrying the labels,
and then remove the indices of labels.

\section{Proof of Theorem \ref{maintheorem}}

We construct all polygonal presentations with $k=3$ and $n=1$ and for
which the graph $G_1$ is a generalized 4-gon. The 23 torsion free
groups were listed in \cite{KV}. Here we give the groups with
torsion. Our strategy is to go through all possible incidence
tableaus for $G_1$ and determine if they can be interpreted as
triangle presentations. 

Let $P$ be the set of black vertices and $L$ be the set of 
white vertices in $G_1$. We denote the elements of $P$ by $x_i$ and the
elements of $Q$ by $y_i$, $i=1,2, \ldots, 15$. In all cases we define
the basic bijection $\lambda:P\rightarrow Q$ by $\lambda(x_i)=y_ i$.

By \cite{Tits}, the smallest thick generalized 4-gon can be presented
in the following way: its ``points'' are pairs $(i,j)$, where
$i,j=1,...,6$, $i \neq j$ and ``lines'' are triples
$(i_1,j_1),(i_2,j_2),(i_3,j_3)$ of those pairs, where $i_1,
i_2,i_3,j_1,j_2$ and $j_3$ are all different. Therefore, we build a
tableau as follows: For each row we take three pairs
$(i_1,j_1),(i_2,j_2)$, and $(i_3,j_3)$, where $i_1, i_2, i_3, j_1, j_2$
and $j_3$ are all different and in $1,2, \ldots,6$. These are our
points: $x_1=(1,2)$, $x_2=(1,3)$,..., $x_{15}=(5,6)$.

\begin{table}[ht]
\centering
\caption{Table of points for incidence tableau}
\label{tabdata}
  \begin{tabular}{ | c c c | }
    \hline
$x_1$ &$x_{10}$ &$x_{15}$\\
$x_1$ &$x_{11}$ &$x_{14}$\\
$x_1$ &$x_{12}$ &$x_{13}$\\
$x_2$ &$x_7$ &$x_{15}$\\
$x_2$ &$x_8$ &$x_{14}$\\
$x_2$ &$x_9$ &$x_{13}$\\
$x_3$ &$x_6$ &$x_{15}$\\
$x_3$ &$x_8$ &$x_{12}$\\
$x_3$ &$x_9$ &$x_{11}$\\
$x_4$ &$x_6$ &$x_{14}$\\
$x_4$ &$x_7$ &$x_{12}$\\
$x_4$ &$x_9$ &$x_{10}$\\
$x_5$ &$x_6$ &$x_{13}$\\
$x_5$ &$x_7$ &$x_{11}$\\
$x_5$ &$x_8$ &$x_{10}$\\
      \hline
  \end{tabular}
\end{table}

Next we label the rows in Table \ref{tabdata} by $y_1, \ldots, y_{15}$
in such a way that the result is an incidence tableau that gives a
triangle presentation with the basic bijection $\lambda$.  To obtain
groups with torsion, we demand that at least one of the triangles is
of the form $(x_i,x_i,x_i)$. For example, labeling the rows from top
to bottom by $y_1$, $y_2$, $y_6$, $y_5$, $y_{14}$, $y_{10}$, $y_7$,
$y_8$, $y_{12}$, $y_3$, $y_4$, $y_9$, $y_{15}$, $y_{13}$, and $y_{11}$
gives rise to the presentation $T_{24}$ with the following 17
triangles: $(x_1, x_1, x_1)$, $(x_{10}, x_2, x_1)$, $(x_{15}, x_6,
x_1)$, $(x_{11}, x_5, x_2)$, $(x_{14}, x_{14}, x_2)$, $(x_4, x_7,
x_3)$ $(x_6, x_{12}, x_3)$, $(x_{14}, x_8, x_3)$, $(x_4, x_4, x_4)$,
$(x_{12}, x_9, x_4)$, $(x_7, x_{15}, x_5)$, $(x_{15}, x_{13}, x_5)$,
$(x_{13}, x_7, x_6)$, $(x_8, x_8, x_8)$, $(x_{12}, x_{11}, x_8)$,
$(x_9, x_{10}, x_9)$, and $(x_{13}, x_{11}, x_{10})$. 

The labeling of rows in Table \ref{tabdata} defines the triangles
uniquely: since the last row $x_5$, $x_8$, $x_{10}$ has label $y_{11}$
we know that there are triangles $(x_{11}, x_5,x_a)$,
$(x_{11},x_8,x_b)$ and $(x_{11},x_{10},x_c)$ for some points $x_a$,
$x_b$ and $x_c$. For the first of these triangles the missing point is
$x_a=x_2$, since the line $y_5$ has points $x_2$, $x_7$ and $x_{15}$
and from the lines with those respective numbers only $y_2$ has the
point $x_{11}$. That is, line $y_{11}$ has point $x_5$, line $y_5$ has
point $x_2$ and line $y_2$ has point $x_{11}$, and this gives the
triangle $(x_{11}, x_5,x_2)$. Similarly, we must have $x_b=x_{12}$ and
$x_c=x_{13}$. Going through all the rows we get the triangles for this
presentation. The number of the triangles in each presentation is
either 17 or 19, depending of whether there is 3 or 6 triangles of the
from $(x_i,x_i,x_i)$.

The presentations are searched by a computer program. The program is
written in Fortran in order to keep it fast and simple.  It goes
through all 15! ways to label the rows of the given tableau, and
decides, which of these give an incidence tableau of a triangle
presentation with torsion. The program outputs one representative of each
equivalence class of triangle presentations. We obtain in this way 7159
different equivalence classes of presentations.

For a polygonal presentation $T$, take $N$ ($N=17$ or $19$) oriented
regular hyperbolic triangles with angles $\pi /4$, write words from
the presentation on their boundaries and glue together sides with the
same letters, respecting orientation.  The result is a hyperbolic
polyhedron with one vertex and $N$ triangular faces, and its universal
covering is a triangular hyperbolic building. We can draw the link,
which is a generalized 4-gon, for any of these buildings: for every
triple $(x_i,x_j,x_k)$ the points $y_i$ and $x_j$, as well as $y_j$
and $x_k$ and $y_k$ and $x_i$ are incident in it. The fundamental
group $\Gamma$ of the polyhedron acts simply transitively on vertices
of the building. The group $\Gamma_i$ has 15 generators and $N$
relations, which come naturally from the polygonal presentation $T$.

To distinguish groups $\Gamma_i$, $i=1,\ldots, 7159$ it is sufficient
to distinguish the isometry classes of polyhedra, according to the
Mostow-type rigidity for hyperbolic buildings which was shown, for
example, in \cite{Xie}.

Therefore, we consider dual graphs of index 3 subgroups in order to
see which of these presentations give rise to isometric
polyhedra. First we calculate the index 3 subgroups: we substitute
each triple of the form $(x_i,x_i,x_i)$ in the presentation by
$(x_i^1,x_i^2,x_i^3)$, each $(x_i,x_j,x_k)$ by three triplets
$(x_i^1,x_j^2,x_k^3)$, $(x_i^2,x_j^3,x_k^1)$ and
$(x_i^3,x_j^1,x_k^2)$, and each $(x_i,x_j,x_j)$ similarly by three
triplets $(x_i^1,x_j^2,x_j^3)$, $(x_i^2,x_j^3,x_j^1)$ and
$(x_i^3,x_j^1,x_j^2)$. We then have 45 triangles, which represent the
generators of the index 3 subgroup of $\Gamma$.

We next construct the dual graph for each of these as follows: we take 90
vertices such that first 45 of them (numbered 1-45) correspond to the
edges of the triangles and the second 45 edges (numbered 46-90)
correspond to the faces of the triangles. We add an edge between
vertices $i$ (from 1-45) and $j$ (from 46-90), if edge $i$ was on the
boundary of the face $j$ in a triangle. Thus we obtain trivalent
bipartite graphs with 90 vertices.

With the help of the Computational Algebra System Magma we compare the
dual graphs of the index 3 subgroups and we find that most of
them are isomorphic with some other graph: there are only 168
non-isomorphic dual graphs. Thus we have 168 triangle presentations
which give rise to non-isometric polyhedra.  We then compute links of
order two in buildings defined by our 168 torsion groups and 23
torsion free groups from \cite{KV}.  There are only two non-isomorphic
links of order two and in this case, they are
complete invariants of buildings.

The 168 triangle presentations are listed in Appendix 1 together with
the number  (1) or (2) indicating the type of building.

This completes the proof of Theorem \ref{maintheorem}.

\section{Construction of polyhedra with m-gonal faces using torsion groups}

In \cite{KV} the authors described how to construct buildings with
m-gonal faces, for arbitrary n, starting from torsion free groups
acting on triangular buildings with the smallest possible link.  We
modify this construction to allow torsion groups and an arbitrary
generalized polygon as the link of each vertex.

Given generalized polygon $G$ we shall denote by $G'$ the graph
arising by calling black resp. white vertices of $G$ black resp.
white vertices of $G'$.

Consider a bipartite graph $G$ with black vertices $P=\{x_1,\ldots, ,
x_{k}\}$ and white vertices $L=\{y_1, \ldots, y_{k}\}$ and a subset
$K\subset P\times P\times P$ that defines the triangles. Starting from
this triangular presentation $K$, we construct a polyhedron, whose
faces are $m$-gons and whose $m$-vertices have links $G$ or $G'$.

Let $w=z_1\ldots z_m$ be a word of length $m$ in three letters $a$,
$b$ and $c$. Assume that $z_1=a$, $z_2=b$ and $z_3=c$ and that $w$
does not contain proper powers of the letters $a$, $b$ and $c$, that
is, $z_m\neq a$ and $z_t\neq z_{t+1}$ for all $t=1,\ldots,m-1$.

For each of the triples $(x_i,x_j,x_k)$ in $K$ we take three triples
$(x_i^1,x_j^2,x_k^3)$, $(x_i^2,x_j^3,x_k^1)$ , $(x_i^3,x_j^1,x_k^2)$
if at least two of $x_i,x_j,x_k$ are different, and just one
$(x_i^1,x_i^2,x_i^3)$, if $i=j=k$.  The triples are cyclic, so we can
write them as $(x_i^1,x_j^2,x_k^3)$, $(x_k^1,x_i^2,x_j^3)$ and
$(x_j^1,x_k^2,x_i^3)$. By glueing together triangles with these words
on the boundary, we obtain a polyhedron with triangle faces and 3
vertices, each of them with the graph $G$ as the link of each vertex.

We construct $m$-tuples, one corresponding to each of these new
triples: for triple $(x_\alpha^1,x_\beta^2,x_\gamma^3)$ we define an
$m$-tuple, which corresponds a word $w$ with $a=x_\alpha^1$,
$b=x_\beta^2$ and $c=x_\gamma^3$. We have $m$-tuples whose coordinates start with one of the triples, and then continue with
$m-3$ letters in some order defined by the word $w$ in the
letters $a$, $b$, and $c$.

If we glue the $m$-gons with these words on the boundary together by
their sides labelled with same letters, respecting orientation, we
obtain a polyhedron with $m$-gonal faces and $m$ vertices, which all
have the link $G$ or $G'$. The type of the link can be seen from the
letters of the edges meeting at that
vertex. Set \[Sign(ab)=Sign(bc)=Sign(ca)=1\]
and \[Sign(ba)=Sign(cb)=Sign(ac)=-1.\] Then for vertex $t=1,\ldots,
m-1$ the group $G_t$ of the link is $G$ if $Sign(z_t,z_{t+1})=1$ and
$G'$ if $Sign(z_t,z_{t+1})=-1$. For the last vertex we have $G_m= G$
if $Sign(z_m,a)=1$ and $G'$ if $Sign(z_m,a)=-1$.

We denote the set of $m$-tuples by $T_m$. Thus we have the following

\begin{theorem}
  The above constructed subset $T_m\subset P\times \cdots \times P$ is
  a polygonal presentation. It defines a polyhedron $X$ whose faces
  are $m$-gons and whose $m$ vertices have links $G$ or $G'$.
\end{theorem}

\section*{Acknowledgements} 

The authors would like to thank Donald Cartwright and Tim Steger for
many helpful discussions and an earlier (unpublished) collaboration
with the first author. We are grateful to Fr\'{e}d\'{e}ric Haglund for
his interest and many illuminating discussions. Some of this work was
completed at the IH\'ES in the Spring of 2010. The first and third authors
would like to express their gratitude to IH\'ES for the Institute's
hospitality.

The first author was supported in part by NSF grant DMS-1101282. The
second author was supported by Emil Aaltonen Foundation, and wishes to
thank \'{E}cole Normale Sup\'{e}rieure for providing excellent working
conditions during Spring 2011. The last author was supported by EPSRC
funded project EP/F014945/1.

\section*{Appendix 1}

The torsion free cases $T_1, \ldots,T_{23}$ have been listed in
\cite{KV}, and thus for those we only denote here in Table \ref{tabT}
whether the link of order 2 is isomorphic to that of $T_1$ (case 1) or
$T_2$ (case 2). Then in Table \ref{tabx} we list the labelings of the
rows of Table \ref{tabdata} which give rise to the triangle
presentations with torsion, denoted by $T_{24}, \ldots, T_{191}$.
After the name of the presentation in Table \ref{tabx} there is (1)
resp. (2) if the resulting building is isomorphic with that of $T_1$
resp. $T_2$.

\begin{table}[ht]
\centering
\caption{List of torsion free cases by isomorphism type of the 2-link}
\label{tabT}
  \begin{tabular}{ | l | l | }
    \hline
    (1) & $T_1$, $T_4$, $T_5$, $T_8$, $T_9$, $T_{11}$, $T_{13}$, $T_{15}$, $T_{17}$, $T_{18}$, $T_{19}$, $T_{23}$\\ \hline
    (2) & $T_2$, $T_3$, $T_6$, $T_7$, $T_{10}$, $T_{12}$, $T_{14}$, $T_{16}$, $T_{20}$, $T_{21}$, $T_{22}$ \\ 
      \hline
  \end{tabular}
\end{table}

\newpage

\begin{table}[ht!]
\begin{center}
\caption{List of labelings giving triangle presentations with torsion}
\label{tabx}
  \begin{tabular}{ | l | l | l |}
    \hline
$T_{24}$&(1)&
$y_{ 1}$, $y_{ 2}$, $y_{ 6}$, $y_{ 5}$, $y_{14}$, $y_{10}$, $y_{ 7}$, $y_{ 8}$, $y_{12}$, $y_{ 3}$, $y_{ 4}$, $y_{ 9}$, $y_{15}$, $y_{13}$, $y_{11}$\\
$T_{25}$&(2)&
$y_{ 1}$, $y_{ 2}$, $y_{ 9}$, $y_{10}$, $y_{ 3}$, $y_{ 5}$, $y_{14}$, $y_{ 8}$, $y_{15}$, $y_{13}$, $y_{ 7}$, $y_{ 4}$, $y_{ 6}$, $y_{12}$, $y_{11}$\\
$T_{26}$&(2)&
$y_{ 1}$, $y_{ 2}$, $y_{ 6}$, $y_{ 3}$, $y_{15}$, $y_{ 7}$, $y_{ 8}$, $y_{14}$, $y_{ 9}$, $y_{10}$, $y_{11}$, $y_{12}$, $y_{ 4}$, $y_{ 5}$, $y_{13}$\\
$T_{27}$&(1)&
$y_{ 1}$, $y_{ 2}$, $y_{14}$, $y_{13}$, $y_{10}$, $y_{ 4}$, $y_{ 3}$, $y_{ 6}$, $y_{ 8}$, $y_{15}$, $y_{11}$, $y_{ 9}$, $y_{ 5}$, $y_{ 7}$, $y_{12}$\\
$T_{28}$&(2)&
$y_{ 1}$, $y_{ 2}$, $y_{ 3}$, $y_{ 5}$, $y_{14}$, $y_{15}$, $y_{ 8}$, $y_{ 9}$, $y_{10}$, $y_{ 6}$, $y_{ 7}$, $y_{12}$, $y_{11}$, $y_{ 4}$, $y_{13}$\\
$T_{29}$&(1)&
$y_{ 1}$, $y_{ 4}$, $y_{12}$, $y_{ 3}$, $y_{ 9}$, $y_{ 6}$, $y_{11}$, $y_{ 5}$, $y_{ 7}$, $y_{14}$, $y_{15}$, $y_{10}$, $y_{13}$, $y_{ 2}$, $y_{ 8}$\\
$T_{30}$&(2)&
$y_{ 1}$, $y_{ 2}$, $y_{ 4}$, $y_{ 9}$, $y_{ 6}$, $y_{10}$, $y_{11}$, $y_{13}$, $y_{14}$, $y_{ 3}$, $y_{15}$, $y_{ 7}$, $y_{ 5}$, $y_{12}$, $y_{ 8}$\\
$T_{31}$&(2)&
$y_{ 2}$, $y_{ 6}$, $y_{13}$, $y_{14}$, $y_{ 3}$, $y_{ 8}$, $y_{ 7}$, $y_{ 1}$, $y_{ 9}$, $y_{15}$, $y_{ 4}$, $y_{ 5}$, $y_{12}$, $y_{11}$, $y_{10}$\\
$T_{32}$&(2)&
$y_{ 1}$, $y_{ 2}$, $y_{ 6}$, $y_{ 8}$, $y_{10}$, $y_{ 4}$, $y_{15}$, $y_{14}$, $y_{13}$, $y_{11}$, $y_{ 9}$, $y_{ 7}$, $y_{ 5}$, $y_{12}$, $y_{ 3}$\\
$T_{33}$&(1)&
$y_{ 1}$, $y_{ 2}$, $y_{ 7}$, $y_{ 3}$, $y_{14}$, $y_{15}$, $y_{ 9}$, $y_{11}$, $y_{13}$, $y_{ 5}$, $y_{ 6}$, $y_{ 4}$, $y_{12}$, $y_{10}$, $y_{ 8}$\\
$T_{34}$&(2)&
$y_{ 2}$, $y_{ 7}$, $y_{ 8}$, $y_{ 4}$, $y_{ 6}$, $y_{ 5}$, $y_{ 1}$, $y_{ 3}$, $y_{ 9}$, $y_{14}$, $y_{15}$, $y_{12}$, $y_{13}$, $y_{11}$, $y_{10}$\\
$T_{35}$&(2)&
$y_{ 4}$, $y_{ 9}$, $y_{12}$, $y_{11}$, $y_{ 8}$, $y_{15}$, $y_{ 1}$, $y_{ 3}$, $y_{ 2}$, $y_{ 6}$, $y_{13}$, $y_{10}$, $y_{ 5}$, $y_{ 7}$, $y_{14}$\\
$T_{36}$&(1)&
$y_{ 1}$, $y_{10}$, $y_{11}$, $y_{ 2}$, $y_{ 8}$, $y_{ 5}$, $y_{14}$, $y_{ 7}$, $y_{ 6}$, $y_{ 9}$, $y_{12}$, $y_{ 3}$, $y_{13}$, $y_{ 4}$, $y_{15}$\\
$T_{37}$&(2)&
$y_{ 1}$, $y_{ 2}$, $y_{ 3}$, $y_{15}$, $y_{ 9}$, $y_{ 4}$, $y_{ 7}$, $y_{12}$, $y_{10}$, $y_{11}$, $y_{ 5}$, $y_{14}$, $y_{ 6}$, $y_{13}$, $y_{ 8}$\\
$T_{38}$&(1)&
$y_{ 1}$, $y_{ 2}$, $y_{ 4}$, $y_{15}$, $y_{13}$, $y_{ 6}$, $y_{11}$, $y_{ 5}$, $y_{12}$, $y_{ 8}$, $y_{ 3}$, $y_{10}$, $y_{14}$, $y_{ 7}$, $y_{ 9}$\\
$T_{39}$&(1)&
$y_{ 1}$, $y_{ 2}$, $y_{ 8}$, $y_{ 4}$, $y_{ 6}$, $y_{10}$, $y_{13}$, $y_{ 7}$, $y_{ 3}$, $y_{14}$, $y_{12}$, $y_{ 5}$, $y_{11}$, $y_{ 9}$, $y_{15}$\\
$T_{40}$&(2)&
$y_{ 1}$, $y_{ 2}$, $y_{ 4}$, $y_{15}$, $y_{14}$, $y_{ 3}$, $y_{ 5}$, $y_{11}$, $y_{ 9}$, $y_{10}$, $y_{12}$, $y_{ 8}$, $y_{ 6}$, $y_{ 7}$, $y_{13}$\\
$T_{41}$&(2)&
$y_{ 1}$, $y_{ 2}$, $y_{ 6}$, $y_{ 4}$, $y_{14}$, $y_{10}$, $y_{ 3}$, $y_{12}$, $y_{13}$, $y_{11}$, $y_{ 7}$, $y_{ 5}$, $y_{15}$, $y_{ 8}$, $y_{ 9}$\\
$T_{42}$&(1)&
$y_{ 1}$, $y_{ 2}$, $y_{ 6}$, $y_{10}$, $y_{13}$, $y_{ 3}$, $y_{ 5}$, $y_{ 8}$, $y_{11}$, $y_{15}$, $y_{ 7}$, $y_{ 4}$, $y_{14}$, $y_{ 9}$, $y_{12}$\\
$T_{43}$&(1)&
$y_{ 1}$, $y_{ 2}$, $y_{ 8}$, $y_{13}$, $y_{15}$, $y_{ 3}$, $y_{ 6}$, $y_{11}$, $y_{ 7}$, $y_{ 9}$, $y_{ 4}$, $y_{ 5}$, $y_{14}$, $y_{12}$, $y_{10}$\\
$T_{44}$&(1)&
$y_{ 1}$, $y_{ 2}$, $y_{ 4}$, $y_{ 5}$, $y_{ 3}$, $y_{15}$, $y_{ 6}$, $y_{11}$, $y_{ 7}$, $y_{ 8}$, $y_{12}$, $y_{10}$, $y_{ 9}$, $y_{14}$, $y_{13}$\\
$T_{45}$&(2)&
$y_{ 1}$, $y_{ 2}$, $y_{ 4}$, $y_{ 8}$, $y_{11}$, $y_{10}$, $y_{13}$, $y_{ 9}$, $y_{ 7}$, $y_{ 6}$, $y_{12}$, $y_{15}$, $y_{ 5}$, $y_{14}$, $y_{ 3}$\\
$T_{46}$&(1)&
$y_{ 1}$, $y_{ 2}$, $y_{ 6}$, $y_{ 3}$, $y_{ 8}$, $y_{10}$, $y_{15}$, $y_{ 5}$, $y_{11}$, $y_{ 9}$, $y_{13}$, $y_{12}$, $y_{ 7}$, $y_{ 4}$, $y_{14}$\\
$T_{47}$&(2)&
$y_{ 1}$, $y_{ 2}$, $y_{ 8}$, $y_{15}$, $y_{ 4}$, $y_{ 6}$, $y_{ 7}$, $y_{ 5}$, $y_{ 3}$, $y_{14}$, $y_{11}$, $y_{12}$, $y_{13}$, $y_{ 9}$, $y_{10}$\\
$T_{48}$&(2)&
$y_{ 1}$, $y_{ 2}$, $y_{ 8}$, $y_{15}$, $y_{ 6}$, $y_{ 7}$, $y_{13}$, $y_{ 4}$, $y_{ 9}$, $y_{12}$, $y_{14}$, $y_{ 5}$, $y_{11}$, $y_{ 3}$, $y_{10}$\\
$T_{49}$&(2)&
$y_{ 1}$, $y_{ 2}$, $y_{ 3}$, $y_{ 5}$, $y_{ 4}$, $y_{15}$, $y_{ 6}$, $y_{10}$, $y_{ 8}$, $y_{ 9}$, $y_{ 7}$, $y_{11}$, $y_{12}$, $y_{14}$, $y_{13}$\\
$T_{50}$&(2)&
$y_{ 1}$, $y_{ 2}$, $y_{ 3}$, $y_{ 6}$, $y_{12}$, $y_{10}$, $y_{ 8}$, $y_{11}$, $y_{15}$, $y_{14}$, $y_{ 4}$, $y_{ 5}$, $y_{ 9}$, $y_{13}$, $y_{ 7}$\\
$T_{51}$&(2)&
$y_{ 1}$, $y_{ 2}$, $y_{ 4}$, $y_{12}$, $y_{ 3}$, $y_{15}$, $y_{ 5}$, $y_{ 9}$, $y_{14}$, $y_{ 7}$, $y_{11}$, $y_{10}$, $y_{ 6}$, $y_{13}$, $y_{ 8}$\\
$T_{52}$&(2)&
$y_{ 1}$, $y_{ 2}$, $y_{ 6}$, $y_{ 3}$, $y_{14}$, $y_{15}$, $y_{10}$, $y_{13}$, $y_{11}$, $y_{12}$, $y_{ 8}$, $y_{ 7}$, $y_{ 4}$, $y_{ 5}$, $y_{ 9}$\\
$T_{53}$&(1)&
$y_{ 1}$, $y_{ 2}$, $y_{ 6}$, $y_{ 9}$, $y_{12}$, $y_{10}$, $y_{15}$, $y_{ 4}$, $y_{ 3}$, $y_{ 5}$, $y_{14}$, $y_{11}$, $y_{ 8}$, $y_{13}$, $y_{ 7}$\\
$T_{54}$&(2)&
$y_{ 1}$, $y_{ 2}$, $y_{12}$, $y_{15}$, $y_{ 9}$, $y_{11}$, $y_{ 7}$, $y_{ 8}$, $y_{ 5}$, $y_{ 6}$, $y_{10}$, $y_{ 3}$, $y_{ 4}$, $y_{14}$, $y_{13}$\\
$T_{55}$&(2)&
$y_{ 1}$, $y_{ 2}$, $y_{13}$, $y_{12}$, $y_{15}$, $y_{ 7}$, $y_{ 9}$, $y_{11}$, $y_{ 4}$, $y_{ 3}$, $y_{14}$, $y_{ 8}$, $y_{10}$, $y_{ 6}$, $y_{ 5}$\\
$T_{56}$&(1)&
$y_{ 2}$, $y_{ 9}$, $y_{15}$, $y_{12}$, $y_{ 8}$, $y_{ 5}$, $y_{ 6}$, $y_{ 1}$, $y_{ 3}$, $y_{11}$, $y_{ 7}$, $y_{ 4}$, $y_{13}$, $y_{10}$, $y_{14}$\\
$T_{57}$&(2)&
$y_{ 1}$, $y_{ 2}$, $y_{ 4}$, $y_{ 5}$, $y_{ 3}$, $y_{15}$, $y_{11}$, $y_{ 7}$, $y_{ 9}$, $y_{ 6}$, $y_{10}$, $y_{ 8}$, $y_{12}$, $y_{14}$, $y_{13}$\\
$T_{58}$&(2)&
$y_{ 1}$, $y_{ 2}$, $y_{ 6}$, $y_{ 5}$, $y_{14}$, $y_{15}$, $y_{10}$, $y_{ 3}$, $y_{ 7}$, $y_{ 4}$, $y_{ 8}$, $y_{12}$, $y_{11}$, $y_{13}$, $y_{ 9}$\\
$T_{59}$&(1)&
$y_{ 1}$, $y_{ 2}$, $y_{ 4}$, $y_{ 8}$, $y_{11}$, $y_{15}$, $y_{ 6}$, $y_{ 9}$, $y_{ 7}$, $y_{10}$, $y_{12}$, $y_{ 5}$, $y_{13}$, $y_{14}$, $y_{ 3}$\\
$T_{60}$&(1)&
$y_{ 1}$, $y_{ 2}$, $y_{12}$, $y_{15}$, $y_{ 9}$, $y_{ 4}$, $y_{ 7}$, $y_{10}$, $y_{13}$, $y_{ 3}$, $y_{11}$, $y_{14}$, $y_{ 8}$, $y_{ 6}$, $y_{ 5}$\\
$T_{61}$&(2)&
$y_{ 1}$, $y_{ 2}$, $y_{12}$, $y_{ 5}$, $y_{14}$, $y_{15}$, $y_{ 8}$, $y_{13}$, $y_{ 3}$, $y_{ 6}$, $y_{10}$, $y_{ 7}$, $y_{11}$, $y_{ 4}$, $y_{ 9}$\\
$T_{62}$&(1)&
$y_{ 1}$, $y_{ 2}$, $y_{ 5}$, $y_{ 6}$, $y_{ 4}$, $y_{15}$, $y_{ 9}$, $y_{ 8}$, $y_{13}$, $y_{14}$, $y_{12}$, $y_{11}$, $y_{10}$, $y_{ 3}$, $y_{ 7}$\\
$T_{63}$&(1)&
$y_{ 1}$, $y_{ 3}$, $y_{ 4}$, $y_{ 8}$, $y_{11}$, $y_{ 9}$, $y_{15}$, $y_{ 6}$, $y_{ 2}$, $y_{10}$, $y_{ 7}$, $y_{ 5}$, $y_{14}$, $y_{13}$, $y_{12}$\\
$T_{64}$&(2)&
$y_{ 1}$, $y_{ 2}$, $y_{ 6}$, $y_{ 7}$, $y_{15}$, $y_{ 3}$, $y_{11}$, $y_{ 9}$, $y_{12}$, $y_{10}$, $y_{ 8}$, $y_{13}$, $y_{ 4}$, $y_{14}$, $y_{ 5}$\\
$T_{65}$&(2)&
$y_{ 1}$, $y_{ 2}$, $y_{ 3}$, $y_{ 8}$, $y_{ 6}$, $y_{10}$, $y_{15}$, $y_{ 9}$, $y_{ 7}$, $y_{14}$, $y_{ 5}$, $y_{12}$, $y_{ 4}$, $y_{13}$, $y_{11}$\\
\hline
\end{tabular}
\end{center}
\end{table}

\begin{center}
    \begin{tabular}{ | l | l | l |}
    \hline 
$T_{66}$&(2)&
$y_{ 1}$, $y_{ 2}$, $y_{ 3}$, $y_{ 4}$, $y_{14}$, $y_{10}$, $y_{15}$, $y_{ 5}$, $y_{ 7}$, $y_{ 6}$, $y_{ 9}$, $y_{11}$, $y_{ 8}$, $y_{13}$, $y_{12}$\\
$T_{67}$&(1)&
$y_{ 1}$, $y_{ 2}$, $y_{ 7}$, $y_{10}$, $y_{14}$, $y_{ 8}$, $y_{ 3}$, $y_{11}$, $y_{ 4}$, $y_{13}$, $y_{ 6}$, $y_{12}$, $y_{ 9}$, $y_{15}$, $y_{ 5}$\\
$T_{68}$&(1)&
$y_{ 1}$, $y_{ 2}$, $y_{ 7}$, $y_{15}$, $y_{ 5}$, $y_{ 4}$, $y_{13}$, $y_{12}$, $y_{ 8}$, $y_{11}$, $y_{ 3}$, $y_{ 6}$, $y_{ 9}$, $y_{10}$, $y_{14}$\\
$T_{69}$&(2)&
$y_{ 1}$, $y_{ 2}$, $y_{ 3}$, $y_{10}$, $y_{14}$, $y_{ 8}$, $y_{ 5}$, $y_{ 4}$, $y_{15}$, $y_{13}$, $y_{ 6}$, $y_{ 9}$, $y_{12}$, $y_{ 7}$, $y_{11}$\\
$T_{70}$&(2)&
$y_{ 1}$, $y_{ 2}$, $y_{ 4}$, $y_{ 6}$, $y_{ 9}$, $y_{15}$, $y_{14}$, $y_{ 8}$, $y_{12}$, $y_{ 5}$, $y_{10}$, $y_{11}$, $y_{13}$, $y_{ 7}$, $y_{ 3}$\\
$T_{71}$&(2)&
$y_{ 1}$, $y_{ 2}$, $y_{11}$, $y_{15}$, $y_{12}$, $y_{ 6}$, $y_{ 7}$, $y_{ 4}$, $y_{ 8}$, $y_{14}$, $y_{ 9}$, $y_{ 3}$, $y_{ 5}$, $y_{10}$, $y_{13}$\\
$T_{72}$&(1)&
$y_{ 1}$, $y_{ 2}$, $y_{12}$, $y_{ 5}$, $y_{ 9}$, $y_{10}$, $y_{ 7}$, $y_{15}$, $y_{11}$, $y_{13}$, $y_{ 4}$, $y_{ 6}$, $y_{ 8}$, $y_{ 3}$, $y_{14}$\\
$T_{73}$&(2)&
$y_{ 2}$, $y_{ 6}$, $y_{13}$, $y_{15}$, $y_{ 3}$, $y_{ 8}$, $y_{ 4}$, $y_{ 1}$, $y_{ 9}$, $y_{14}$, $y_{11}$, $y_{ 7}$, $y_{12}$, $y_{ 5}$, $y_{10}$\\
$T_{74}$&(2)&
$y_{ 1}$, $y_{ 2}$, $y_{ 4}$, $y_{ 9}$, $y_{ 8}$, $y_{10}$, $y_{ 7}$, $y_{ 3}$, $y_{14}$, $y_{ 5}$, $y_{12}$, $y_{15}$, $y_{ 6}$, $y_{13}$, $y_{11}$\\
$T_{75}$&(1)&
$y_{ 1}$, $y_{10}$, $y_{11}$, $y_{ 5}$, $y_{14}$, $y_{13}$, $y_{ 6}$, $y_{ 7}$, $y_{ 2}$, $y_{ 3}$, $y_{12}$, $y_{15}$, $y_{ 9}$, $y_{ 4}$, $y_{ 8}$\\
$T_{76}$&(1)&
$y_{ 1}$, $y_{ 2}$, $y_{12}$, $y_{10}$, $y_{14}$, $y_{ 5}$, $y_{ 9}$, $y_{13}$, $y_{ 7}$, $y_{ 6}$, $y_{15}$, $y_{ 4}$, $y_{ 8}$, $y_{11}$, $y_{ 3}$\\
$T_{77}$&(2)&
$y_{ 1}$, $y_{ 2}$, $y_{ 4}$, $y_{ 7}$, $y_{10}$, $y_{ 9}$, $y_{ 5}$, $y_{ 3}$, $y_{13}$, $y_{15}$, $y_{ 6}$, $y_{11}$, $y_{12}$, $y_{14}$, $y_{ 8}$\\
$T_{78}$&(1)&
$y_{ 1}$, $y_{ 2}$, $y_{ 5}$, $y_{ 8}$, $y_{12}$, $y_{10}$, $y_{ 4}$, $y_{ 7}$, $y_{ 3}$, $y_{13}$, $y_{11}$, $y_{ 9}$, $y_{ 6}$, $y_{15}$, $y_{14}$\\
$T_{79}$&(1)&
$y_{ 1}$, $y_{ 2}$, $y_{ 4}$, $y_{ 3}$, $y_{14}$, $y_{15}$, $y_{ 5}$, $y_{11}$, $y_{ 9}$, $y_{10}$, $y_{12}$, $y_{ 8}$, $y_{ 6}$, $y_{ 7}$, $y_{13}$\\
$T_{80}$&(2)&
$y_{ 1}$, $y_{ 2}$, $y_{ 5}$, $y_{12}$, $y_{14}$, $y_{10}$, $y_{13}$, $y_{11}$, $y_{ 6}$, $y_{ 4}$, $y_{ 8}$, $y_{ 9}$, $y_{ 3}$, $y_{ 7}$, $y_{15}$\\
$T_{81}$&(2)&
$y_{ 1}$, $y_{ 2}$, $y_{ 8}$, $y_{ 4}$, $y_{14}$, $y_{15}$, $y_{12}$, $y_{ 6}$, $y_{ 7}$, $y_{13}$, $y_{ 9}$, $y_{11}$, $y_{ 3}$, $y_{ 5}$, $y_{10}$\\
$T_{82}$&(1)&
$y_{ 1}$, $y_{ 2}$, $y_{ 4}$, $y_{ 9}$, $y_{15}$, $y_{11}$, $y_{ 7}$, $y_{12}$, $y_{14}$, $y_{10}$, $y_{ 8}$, $y_{ 3}$, $y_{ 6}$, $y_{ 5}$, $y_{13}$\\
$T_{83}$&(1)&
$y_{ 1}$, $y_{ 2}$, $y_{ 4}$, $y_{13}$, $y_{15}$, $y_{ 8}$, $y_{ 6}$, $y_{ 5}$, $y_{12}$, $y_{ 3}$, $y_{ 7}$, $y_{10}$, $y_{14}$, $y_{ 9}$, $y_{11}$\\
$T_{84}$&(2)&
$y_{ 1}$, $y_{ 2}$, $y_{ 8}$, $y_{ 6}$, $y_{ 4}$, $y_{15}$, $y_{ 9}$, $y_{ 3}$, $y_{13}$, $y_{11}$, $y_{ 5}$, $y_{12}$, $y_{14}$, $y_{ 7}$, $y_{10}$\\
$T_{85}$&(1)&
$y_{ 1}$, $y_{ 2}$, $y_{11}$, $y_{15}$, $y_{14}$, $y_{13}$, $y_{ 7}$, $y_{ 3}$, $y_{10}$, $y_{ 5}$, $y_{ 4}$, $y_{ 6}$, $y_{ 9}$, $y_{ 8}$, $y_{12}$\\
$T_{86}$&(2)&
$y_{ 1}$, $y_{ 6}$, $y_{ 8}$, $y_{11}$, $y_{ 9}$, $y_{13}$, $y_{15}$, $y_{12}$, $y_{ 4}$, $y_{14}$, $y_{ 7}$, $y_{ 5}$, $y_{ 2}$, $y_{ 3}$, $y_{10}$\\
$T_{87}$&(1)&
$y_{ 1}$, $y_{ 2}$, $y_{ 4}$, $y_{ 3}$, $y_{ 6}$, $y_{15}$, $y_{ 5}$, $y_{11}$, $y_{ 9}$, $y_{10}$, $y_{12}$, $y_{ 8}$, $y_{14}$, $y_{ 7}$, $y_{13}$\\
$T_{88}$&(2)&
$y_{ 1}$, $y_{ 2}$, $y_{14}$, $y_{15}$, $y_{ 5}$, $y_{13}$, $y_{ 3}$, $y_{ 4}$, $y_{ 7}$, $y_{10}$, $y_{ 6}$, $y_{ 8}$, $y_{ 9}$, $y_{11}$, $y_{12}$\\
$T_{89}$&(1)&
$y_{ 1}$, $y_{ 2}$, $y_{ 3}$, $y_{10}$, $y_{13}$, $y_{ 5}$, $y_{ 7}$, $y_{ 8}$, $y_{15}$, $y_{ 4}$, $y_{ 9}$, $y_{ 6}$, $y_{14}$, $y_{12}$, $y_{11}$\\
$T_{90}$&(2)&
$y_{ 1}$, $y_{ 2}$, $y_{ 4}$, $y_{ 5}$, $y_{14}$, $y_{15}$, $y_{ 8}$, $y_{12}$, $y_{ 7}$, $y_{ 3}$, $y_{ 6}$, $y_{10}$, $y_{13}$, $y_{11}$, $y_{ 9}$\\
$T_{91}$&(2)&
$y_{ 1}$, $y_{ 2}$, $y_{ 4}$, $y_{ 6}$, $y_{15}$, $y_{12}$, $y_{ 8}$, $y_{14}$, $y_{13}$, $y_{11}$, $y_{ 9}$, $y_{10}$, $y_{ 5}$, $y_{ 3}$, $y_{ 7}$\\
$T_{92}$&(1)&
$y_{ 1}$, $y_{ 2}$, $y_{ 6}$, $y_{10}$, $y_{ 9}$, $y_{ 4}$, $y_{15}$, $y_{12}$, $y_{14}$, $y_{ 8}$, $y_{ 5}$, $y_{11}$, $y_{13}$, $y_{ 7}$, $y_{ 3}$\\
$T_{93}$&(1)&
$y_{ 1}$, $y_{ 2}$, $y_{ 9}$, $y_{14}$, $y_{10}$, $y_{ 7}$, $y_{12}$, $y_{ 4}$, $y_{ 3}$, $y_{11}$, $y_{ 5}$, $y_{15}$, $y_{ 6}$, $y_{13}$, $y_{ 8}$\\
$T_{94}$&(2)&
$y_{ 1}$, $y_{ 2}$, $y_{11}$, $y_{10}$, $y_{ 3}$, $y_{12}$, $y_{14}$, $y_{ 8}$, $y_{15}$, $y_{ 7}$, $y_{ 5}$, $y_{ 4}$, $y_{13}$, $y_{ 6}$, $y_{ 9}$\\
$T_{95}$&(2)&
$y_{ 2}$, $y_{14}$, $y_{13}$, $y_{ 6}$, $y_{ 3}$, $y_{ 4}$, $y_{15}$, $y_{ 8}$, $y_{11}$, $y_{ 5}$, $y_{12}$, $y_{10}$, $y_{ 1}$, $y_{ 9}$, $y_{ 7}$\\
$T_{96}$&(2)&
$y_{ 1}$, $y_{ 3}$, $y_{14}$, $y_{ 9}$, $y_{15}$, $y_{ 2}$, $y_{10}$, $y_{12}$, $y_{11}$, $y_{ 4}$, $y_{ 5}$, $y_{13}$, $y_{ 6}$, $y_{ 7}$, $y_{ 8}$\\
$T_{97}$&(2)&
$y_{ 2}$, $y_{ 3}$, $y_{13}$, $y_{15}$, $y_{ 4}$, $y_{ 8}$, $y_{ 6}$, $y_{12}$, $y_{ 5}$, $y_{14}$, $y_{ 1}$, $y_{ 9}$, $y_{ 7}$, $y_{11}$, $y_{10}$\\
$T_{98}$&(1)&
$y_{ 2}$, $y_{10}$, $y_{15}$, $y_{13}$, $y_{ 9}$, $y_{ 1}$, $y_{ 6}$, $y_{ 8}$, $y_{11}$, $y_{ 4}$, $y_{ 7}$, $y_{ 3}$, $y_{12}$, $y_{14}$, $y_{ 5}$\\
$T_{99}$&(2)&
$y_{ 2}$, $y_{11}$, $y_{15}$, $y_{ 8}$, $y_{ 1}$, $y_{ 9}$, $y_{ 6}$, $y_{12}$, $y_{14}$, $y_{13}$, $y_{ 7}$, $y_{10}$, $y_{ 5}$, $y_{ 3}$, $y_{ 4}$\\
$T_{100}$&(1)&
$y_{ 2}$, $y_{14}$, $y_{15}$, $y_{12}$, $y_{ 7}$, $y_{ 5}$, $y_{ 6}$, $y_{ 8}$, $y_{ 9}$, $y_{ 4}$, $y_{ 1}$, $y_{ 3}$, $y_{13}$, $y_{11}$, $y_{10}$\\
$T_{101}$&(1)&
$y_{ 2}$, $y_{ 4}$, $y_{15}$, $y_{ 8}$, $y_{13}$, $y_{ 1}$, $y_{ 3}$, $y_{12}$, $y_{ 5}$, $y_{14}$, $y_{ 7}$, $y_{ 9}$, $y_{ 6}$, $y_{11}$, $y_{10}$\\
$T_{102}$&(2)&
$y_{ 1}$, $y_{10}$, $y_{11}$, $y_{ 7}$, $y_{ 5}$, $y_{ 9}$, $y_{13}$, $y_{ 2}$, $y_{ 3}$, $y_{14}$, $y_{12}$, $y_{15}$, $y_{ 8}$, $y_{ 4}$, $y_{ 6}$\\
$T_{103}$&(2)&
$y_{ 2}$, $y_{ 4}$, $y_{13}$, $y_{ 7}$, $y_{ 6}$, $y_{ 9}$, $y_{ 8}$, $y_{14}$, $y_{12}$, $y_{ 3}$, $y_{ 5}$, $y_{10}$, $y_{ 1}$, $y_{15}$, $y_{11}$\\
$T_{104}$&(1)&
$y_{ 2}$, $y_{10}$, $y_{12}$, $y_{15}$, $y_{ 1}$, $y_{14}$, $y_{ 6}$, $y_{ 4}$, $y_{11}$, $y_{ 7}$, $y_{ 8}$, $y_{13}$, $y_{ 5}$, $y_{ 9}$, $y_{ 3}$\\
$T_{105}$&(1)&
$y_{ 2}$, $y_{ 6}$, $y_{11}$, $y_{ 8}$, $y_{13}$, $y_{ 7}$, $y_{ 9}$, $y_{12}$, $y_{14}$, $y_{ 5}$, $y_{15}$, $y_{ 4}$, $y_{ 1}$, $y_{ 3}$, $y_{10}$\\
$T_{106}$&(1)&
$y_{ 3}$, $y_{ 6}$, $y_{13}$, $y_{ 9}$, $y_{ 2}$, $y_{ 4}$, $y_{ 7}$, $y_{10}$, $y_{ 8}$, $y_{14}$, $y_{ 1}$, $y_{15}$, $y_{12}$, $y_{11}$, $y_{ 5}$\\
$T_{107}$&(2)&
$y_{ 2}$, $y_{ 3}$, $y_{13}$, $y_{ 8}$, $y_{15}$, $y_{ 7}$, $y_{10}$, $y_{ 9}$, $y_{11}$, $y_{ 6}$, $y_{12}$, $y_{ 4}$, $y_{ 5}$, $y_{ 1}$, $y_{14}$\\
\hline
\end{tabular}

    \begin{tabular}{ | l | l | l |}
    \hline 
$T_{108}$&(1)&
$y_{ 2}$, $y_{ 3}$, $y_{10}$, $y_{13}$, $y_{ 8}$, $y_{12}$, $y_{11}$, $y_{15}$, $y_{ 4}$, $y_{ 6}$, $y_{ 1}$, $y_{ 9}$, $y_{ 5}$, $y_{14}$, $y_{ 7}$\\
$T_{109}$&(1)&
$y_{ 2}$, $y_{ 7}$, $y_{13}$, $y_{ 6}$, $y_{ 3}$, $y_{ 9}$, $y_{15}$, $y_{ 5}$, $y_{12}$, $y_{ 1}$, $y_{ 8}$, $y_{10}$, $y_{14}$, $y_{ 4}$, $y_{11}$\\
$T_{110}$&(1)&
$y_{ 1}$, $y_{ 2}$, $y_{ 4}$, $y_{15}$, $y_{13}$, $y_{12}$, $y_{ 3}$, $y_{ 8}$, $y_{ 7}$, $y_{ 9}$, $y_{14}$, $y_{10}$, $y_{11}$, $y_{ 5}$, $y_{ 6}$\\
$T_{111}$&(1)&
$y_{ 1}$, $y_{ 2}$, $y_{ 3}$, $y_{ 8}$, $y_{ 9}$, $y_{10}$, $y_{15}$, $y_{ 6}$, $y_{11}$, $y_{12}$, $y_{ 5}$, $y_{ 7}$, $y_{ 4}$, $y_{13}$, $y_{14}$\\
$T_{112}$&(2)&
$y_{ 1}$, $y_{ 2}$, $y_{ 3}$, $y_{10}$, $y_{14}$, $y_{ 8}$, $y_{ 6}$, $y_{ 5}$, $y_{15}$, $y_{13}$, $y_{ 7}$, $y_{ 4}$, $y_{ 9}$, $y_{12}$, $y_{11}$\\
$T_{113}$&(2)&
$y_{ 1}$, $y_{ 2}$, $y_{ 4}$, $y_{15}$, $y_{ 6}$, $y_{ 3}$, $y_{ 5}$, $y_{14}$, $y_{12}$, $y_{11}$, $y_{ 9}$, $y_{10}$, $y_{ 8}$, $y_{13}$, $y_{ 7}$\\
$T_{114}$&(2)&
$y_{ 1}$, $y_{ 2}$, $y_{ 4}$, $y_{15}$, $y_{13}$, $y_{ 5}$, $y_{ 9}$, $y_{ 3}$, $y_{ 7}$, $y_{10}$, $y_{12}$, $y_{ 8}$, $y_{14}$, $y_{11}$, $y_{ 6}$\\
$T_{115}$&(1)&
$y_{ 1}$, $y_{ 2}$, $y_{ 4}$, $y_{15}$, $y_{13}$, $y_{11}$, $y_{12}$, $y_{ 7}$, $y_{14}$, $y_{ 8}$, $y_{ 6}$, $y_{10}$, $y_{ 3}$, $y_{ 5}$, $y_{ 9}$\\
$T_{116}$&(1)&
$y_{ 1}$, $y_{ 2}$, $y_{ 4}$, $y_{15}$, $y_{14}$, $y_{ 6}$, $y_{ 8}$, $y_{13}$, $y_{ 9}$, $y_{ 3}$, $y_{10}$, $y_{ 5}$, $y_{11}$, $y_{12}$, $y_{ 7}$\\
$T_{117}$&(1)&
$y_{ 1}$, $y_{ 2}$, $y_{ 6}$, $y_{ 8}$, $y_{13}$, $y_{15}$, $y_{10}$, $y_{ 4}$, $y_{ 3}$, $y_{ 7}$, $y_{ 5}$, $y_{ 9}$, $y_{14}$, $y_{12}$, $y_{11}$\\
$T_{118}$&(1)&
$y_{ 1}$, $y_{ 6}$, $y_{11}$, $y_{ 3}$, $y_{ 9}$, $y_{12}$, $y_{15}$, $y_{ 5}$, $y_{10}$, $y_{ 4}$, $y_{14}$, $y_{ 7}$, $y_{13}$, $y_{ 2}$, $y_{ 8}$\\
$T_{119}$&(2)&
$y_{ 1}$, $y_{ 2}$, $y_{ 4}$, $y_{10}$, $y_{ 9}$, $y_{11}$, $y_{13}$, $y_{ 8}$, $y_{ 3}$, $y_{ 6}$, $y_{15}$, $y_{ 5}$, $y_{ 7}$, $y_{14}$, $y_{12}$\\
$T_{120}$&(1)&
$y_{ 1}$, $y_{ 2}$, $y_{ 6}$, $y_{10}$, $y_{ 3}$, $y_{ 5}$, $y_{ 9}$, $y_{ 8}$, $y_{14}$, $y_{ 7}$, $y_{12}$, $y_{13}$, $y_{15}$, $y_{11}$, $y_{ 4}$\\
$T_{121}$&(2)&
$y_{ 1}$, $y_{ 2}$, $y_{ 8}$, $y_{ 4}$, $y_{ 7}$, $y_{10}$, $y_{ 9}$, $y_{ 3}$, $y_{12}$, $y_{11}$, $y_{14}$, $y_{ 5}$, $y_{ 6}$, $y_{13}$, $y_{15}$\\
$T_{122}$&(2)&
$y_{ 1}$, $y_{ 2}$, $y_{ 3}$, $y_{ 9}$, $y_{15}$, $y_{ 8}$, $y_{10}$, $y_{ 5}$, $y_{ 6}$, $y_{ 7}$, $y_{ 4}$, $y_{14}$, $y_{13}$, $y_{12}$, $y_{11}$\\
$T_{123}$&(1)&
$y_{ 1}$, $y_{ 2}$, $y_{ 6}$, $y_{ 3}$, $y_{10}$, $y_{ 8}$, $y_{ 5}$, $y_{14}$, $y_{11}$, $y_{15}$, $y_{ 9}$, $y_{13}$, $y_{ 4}$, $y_{ 7}$, $y_{12}$\\
$T_{124}$&(2)&
$y_{ 1}$, $y_{ 2}$, $y_{ 8}$, $y_{13}$, $y_{ 7}$, $y_{10}$, $y_{ 3}$, $y_{12}$, $y_{ 5}$, $y_{ 9}$, $y_{ 6}$, $y_{ 4}$, $y_{14}$, $y_{11}$, $y_{15}$\\
$T_{125}$&(2)&
$y_{ 1}$, $y_{ 2}$, $y_{ 4}$, $y_{10}$, $y_{ 8}$, $y_{ 9}$, $y_{12}$, $y_{ 5}$, $y_{ 7}$, $y_{15}$, $y_{ 6}$, $y_{11}$, $y_{ 3}$, $y_{13}$, $y_{14}$\\
$T_{126}$&(2)&
$y_{ 1}$, $y_{ 2}$, $y_{ 6}$, $y_{ 3}$, $y_{10}$, $y_{ 7}$, $y_{15}$, $y_{13}$, $y_{14}$, $y_{ 5}$, $y_{ 8}$, $y_{ 4}$, $y_{12}$, $y_{11}$, $y_{ 9}$\\
$T_{127}$&(2)&
$y_{10}$, $y_{11}$, $y_{13}$, $y_{ 2}$, $y_{ 5}$, $y_{12}$, $y_{14}$, $y_{ 7}$, $y_{ 1}$, $y_{ 8}$, $y_{ 4}$, $y_{ 9}$, $y_{ 3}$, $y_{ 6}$, $y_{15}$\\
$T_{128}$&(1)&
$y_{ 1}$, $y_{ 4}$, $y_{ 7}$, $y_{15}$, $y_{11}$, $y_{ 8}$, $y_{ 6}$, $y_{ 3}$, $y_{13}$, $y_{12}$, $y_{14}$, $y_{10}$, $y_{ 2}$, $y_{ 9}$, $y_{ 5}$\\
$T_{129}$&(1)&
$y_{ 1}$, $y_{ 2}$, $y_{ 6}$, $y_{15}$, $y_{ 5}$, $y_{ 9}$, $y_{ 7}$, $y_{ 8}$, $y_{13}$, $y_{ 3}$, $y_{ 4}$, $y_{11}$, $y_{10}$, $y_{12}$, $y_{14}$\\
$T_{130}$&(1)&
$y_{ 8}$, $y_{14}$, $y_{15}$, $y_{ 7}$, $y_{11}$, $y_{ 6}$, $y_{ 3}$, $y_{12}$, $y_{ 9}$, $y_{ 4}$, $y_{ 1}$, $y_{ 2}$, $y_{ 5}$, $y_{10}$, $y_{13}$\\
$T_{131}$&(2)&
$y_{ 1}$, $y_{ 2}$, $y_{13}$, $y_{15}$, $y_{ 9}$, $y_{ 3}$, $y_{12}$, $y_{14}$, $y_{ 8}$, $y_{ 6}$, $y_{ 7}$, $y_{11}$, $y_{10}$, $y_{ 5}$, $y_{ 4}$\\
$T_{132}$&(1)&
$y_{ 1}$, $y_{ 2}$, $y_{ 3}$, $y_{ 4}$, $y_{10}$, $y_{13}$, $y_{ 7}$, $y_{15}$, $y_{ 9}$, $y_{14}$, $y_{12}$, $y_{ 6}$, $y_{11}$, $y_{ 8}$, $y_{ 5}$\\
$T_{133}$&(2)&
$y_{ 1}$, $y_{ 2}$, $y_{ 4}$, $y_{15}$, $y_{ 3}$, $y_{ 6}$, $y_{11}$, $y_{12}$, $y_{13}$, $y_{10}$, $y_{ 8}$, $y_{ 9}$, $y_{14}$, $y_{ 7}$, $y_{ 5}$\\
$T_{134}$&(1)&
$y_{ 1}$, $y_{ 2}$, $y_{ 8}$, $y_{15}$, $y_{ 6}$, $y_{ 9}$, $y_{ 4}$, $y_{12}$, $y_{11}$, $y_{14}$, $y_{ 7}$, $y_{ 5}$, $y_{13}$, $y_{ 3}$, $y_{10}$\\
$T_{135}$&(2)&
$y_{ 1}$, $y_{ 2}$, $y_{12}$, $y_{15}$, $y_{ 8}$, $y_{ 3}$, $y_{ 7}$, $y_{10}$, $y_{11}$, $y_{ 6}$, $y_{ 4}$, $y_{13}$, $y_{ 5}$, $y_{ 9}$, $y_{14}$\\
$T_{136}$&(1)&
$y_{ 1}$, $y_{ 2}$, $y_{13}$, $y_{ 3}$, $y_{14}$, $y_{15}$, $y_{12}$, $y_{11}$, $y_{ 9}$, $y_{ 4}$, $y_{ 7}$, $y_{ 8}$, $y_{10}$, $y_{ 6}$, $y_{ 5}$\\
$T_{137}$&(1)&
$y_{ 1}$, $y_{ 4}$, $y_{12}$, $y_{ 2}$, $y_{13}$, $y_{ 7}$, $y_{14}$, $y_{ 8}$, $y_{ 9}$, $y_{ 5}$, $y_{15}$, $y_{10}$, $y_{ 6}$, $y_{ 3}$, $y_{11}$\\
$T_{138}$&(1)&
$y_{ 1}$, $y_{ 2}$, $y_{12}$, $y_{10}$, $y_{ 4}$, $y_{ 9}$, $y_{ 7}$, $y_{15}$, $y_{11}$, $y_{14}$, $y_{ 5}$, $y_{ 6}$, $y_{13}$, $y_{ 3}$, $y_{ 8}$\\
$T_{139}$&(1)&
$y_{ 1}$, $y_{ 2}$, $y_{ 3}$, $y_{ 4}$, $y_{10}$, $y_{13}$, $y_{15}$, $y_{ 8}$, $y_{ 9}$, $y_{11}$, $y_{ 7}$, $y_{ 6}$, $y_{14}$, $y_{ 5}$, $y_{12}$\\
$T_{140}$&(2)&
$y_{ 1}$, $y_{ 2}$, $y_{ 4}$, $y_{ 3}$, $y_{14}$, $y_{10}$, $y_{13}$, $y_{ 9}$, $y_{11}$, $y_{ 5}$, $y_{ 6}$, $y_{15}$, $y_{ 7}$, $y_{12}$, $y_{ 8}$\\
$T_{141}$&(2)&
$y_{ 1}$, $y_{ 2}$, $y_{ 4}$, $y_{ 8}$, $y_{14}$, $y_{10}$, $y_{13}$, $y_{ 9}$, $y_{ 7}$, $y_{11}$, $y_{ 6}$, $y_{15}$, $y_{ 5}$, $y_{12}$, $y_{ 3}$\\
$T_{142}$&(1)&
$y_{ 1}$, $y_{ 2}$, $y_{13}$, $y_{ 5}$, $y_{10}$, $y_{15}$, $y_{14}$, $y_{ 8}$, $y_{ 7}$, $y_{ 4}$, $y_{ 6}$, $y_{12}$, $y_{ 9}$, $y_{11}$, $y_{ 3}$\\
$T_{143}$&(2)&
$y_{ 1}$, $y_{ 2}$, $y_{12}$, $y_{ 6}$, $y_{14}$, $y_{15}$, $y_{11}$, $y_{10}$, $y_{ 8}$, $y_{ 9}$, $y_{ 7}$, $y_{ 5}$, $y_{ 4}$, $y_{13}$, $y_{ 3}$\\
$T_{144}$&(2)&
$y_{ 1}$, $y_{ 2}$, $y_{ 4}$, $y_{15}$, $y_{13}$, $y_{ 6}$, $y_{11}$, $y_{ 3}$, $y_{12}$, $y_{10}$, $y_{ 8}$, $y_{ 9}$, $y_{14}$, $y_{ 7}$, $y_{ 5}$\\
$T_{145}$&(1)&
$y_{ 1}$, $y_{ 2}$, $y_{ 6}$, $y_{10}$, $y_{ 5}$, $y_{ 9}$, $y_{ 7}$, $y_{ 8}$, $y_{12}$, $y_{15}$, $y_{13}$, $y_{14}$, $y_{11}$, $y_{ 4}$, $y_{ 3}$\\
$T_{146}$&(2)&
$y_{ 1}$, $y_{ 2}$, $y_{ 9}$, $y_{15}$, $y_{14}$, $y_{ 8}$, $y_{ 5}$, $y_{12}$, $y_{ 6}$, $y_{13}$, $y_{ 3}$, $y_{10}$, $y_{ 7}$, $y_{ 4}$, $y_{11}$\\
$T_{147}$&(2)&
$y_{ 1}$, $y_{ 2}$, $y_{ 4}$, $y_{ 5}$, $y_{12}$, $y_{10}$, $y_{ 3}$, $y_{14}$, $y_{13}$, $y_{15}$, $y_{ 8}$, $y_{ 9}$, $y_{ 6}$, $y_{ 7}$, $y_{11}$\\
$T_{148}$&(1)&
$y_{ 1}$, $y_{ 2}$, $y_{ 5}$, $y_{15}$, $y_{12}$, $y_{ 4}$, $y_{ 6}$, $y_{11}$, $y_{13}$, $y_{ 7}$, $y_{ 9}$, $y_{14}$, $y_{ 8}$, $y_{10}$, $y_{ 3}$\\
$T_{149}$&(2)&
$y_{ 1}$, $y_{10}$, $y_{11}$, $y_{ 3}$, $y_{14}$, $y_{ 7}$, $y_{13}$, $y_{ 8}$, $y_{ 4}$, $y_{12}$, $y_{ 5}$, $y_{15}$, $y_{ 9}$, $y_{ 6}$, $y_{ 2}$\\
\hline
\end{tabular}

    \begin{tabular}{ | l | l | l |}
    \hline 
$T_{150}$&(1)&
$y_{ 5}$, $y_{14}$, $y_{15}$, $y_{ 3}$, $y_{ 8}$, $y_{ 9}$, $y_{11}$, $y_{12}$, $y_{ 7}$, $y_{ 6}$, $y_{ 2}$, $y_{ 4}$, $y_{10}$, $y_{13}$, $y_{ 1}$\\
$T_{151}$&(1)&
$y_{ 6}$, $y_{ 8}$, $y_{10}$, $y_{ 3}$, $y_{14}$, $y_{15}$, $y_{ 4}$, $y_{ 2}$, $y_{12}$, $y_{ 1}$, $y_{ 7}$, $y_{ 9}$, $y_{13}$, $y_{11}$, $y_{ 5}$\\
$T_{152}$&(1)&
$y_{ 4}$, $y_{11}$, $y_{12}$, $y_{15}$, $y_{ 6}$, $y_{13}$, $y_{ 7}$, $y_{ 8}$, $y_{ 9}$, $y_{10}$, $y_{ 2}$, $y_{14}$, $y_{ 3}$, $y_{ 5}$, $y_{ 1}$\\
$T_{153}$&(2)&
$y_{ 6}$, $y_{11}$, $y_{14}$, $y_{10}$, $y_{ 8}$, $y_{ 2}$, $y_{15}$, $y_{ 9}$, $y_{ 5}$, $y_{ 7}$, $y_{12}$, $y_{ 4}$, $y_{ 3}$, $y_{ 1}$, $y_{13}$\\
$T_{154}$&(2)&
$y_{ 2}$, $y_{ 5}$, $y_{11}$, $y_{13}$, $y_{ 8}$, $y_{ 9}$, $y_{12}$, $y_{ 4}$, $y_{15}$, $y_{14}$, $y_{ 3}$, $y_{ 7}$, $y_{10}$, $y_{ 6}$, $y_{ 1}$\\
$T_{155}$&(1)&
$y_{ 1}$, $y_{ 8}$, $y_{12}$, $y_{11}$, $y_{14}$, $y_{13}$, $y_{ 6}$, $y_{ 2}$, $y_{ 5}$, $y_{ 3}$, $y_{15}$, $y_{ 9}$, $y_{ 4}$, $y_{ 7}$, $y_{10}$\\
$T_{156}$&(2)&
$y_{ 1}$, $y_{ 6}$, $y_{12}$, $y_{ 2}$, $y_{ 3}$, $y_{13}$, $y_{10}$, $y_{15}$, $y_{ 5}$, $y_{ 4}$, $y_{ 9}$, $y_{14}$, $y_{ 7}$, $y_{11}$, $y_{ 8}$\\
$T_{157}$&(2)&
$y_{ 8}$, $y_{ 1}$, $y_{ 5}$, $y_{ 3}$, $y_{14}$, $y_{10}$, $y_{15}$, $y_{ 9}$, $y_{ 7}$, $y_{12}$, $y_{ 2}$, $y_{ 4}$, $y_{13}$, $y_{11}$, $y_{ 6}$\\
$T_{158}$&(1)&
$y_{ 6}$, $y_{12}$, $y_{14}$, $y_{15}$, $y_{ 7}$, $y_{ 3}$, $y_{ 8}$, $y_{13}$, $y_{ 9}$, $y_{11}$, $y_{ 5}$, $y_{ 4}$, $y_{ 2}$, $y_{ 1}$, $y_{10}$\\
$T_{159}$&(2)&
$y_{ 6}$, $y_{ 9}$, $y_{12}$, $y_{15}$, $y_{ 7}$, $y_{10}$, $y_{ 2}$, $y_{14}$, $y_{11}$, $y_{13}$, $y_{ 8}$, $y_{ 3}$, $y_{ 5}$, $y_{ 4}$, $y_{ 1}$\\
$T_{160}$&(1)&
$y_{ 2}$, $y_{11}$, $y_{12}$, $y_{ 3}$, $y_{ 5}$, $y_{ 7}$, $y_{15}$, $y_{ 6}$, $y_{ 9}$, $y_{13}$, $y_{ 1}$, $y_{14}$, $y_{10}$, $y_{ 4}$, $y_{ 8}$\\
$T_{161}$&(2)&
$y_{ 2}$, $y_{ 9}$, $y_{13}$, $y_{15}$, $y_{12}$, $y_{ 8}$, $y_{ 6}$, $y_{10}$, $y_{ 5}$, $y_{ 3}$, $y_{ 1}$, $y_{ 4}$, $y_{ 7}$, $y_{11}$, $y_{14}$\\
$T_{162}$&(1)&
$y_{ 2}$, $y_{ 3}$, $y_{13}$, $y_{11}$, $y_{ 8}$, $y_{ 6}$, $y_{15}$, $y_{10}$, $y_{ 5}$, $y_{12}$, $y_{ 9}$, $y_{14}$, $y_{ 7}$, $y_{ 1}$, $y_{ 4}$\\
$T_{163}$&(1)&
$y_{ 1}$, $y_{ 2}$, $y_{ 7}$, $y_{15}$, $y_{ 4}$, $y_{ 9}$, $y_{13}$, $y_{ 8}$, $y_{11}$, $y_{14}$, $y_{12}$, $y_{ 6}$, $y_{ 5}$, $y_{10}$, $y_{ 3}$\\
$T_{164}$&(1)&
$y_{ 2}$, $y_{ 5}$, $y_{10}$, $y_{15}$, $y_{13}$, $y_{ 1}$, $y_{ 6}$, $y_{14}$, $y_{11}$, $y_{ 7}$, $y_{12}$, $y_{ 4}$, $y_{ 3}$, $y_{ 9}$, $y_{ 8}$\\
$T_{165}$&(1)&
$y_{ 2}$, $y_{ 9}$, $y_{10}$, $y_{12}$, $y_{ 6}$, $y_{ 5}$, $y_{ 8}$, $y_{11}$, $y_{ 3}$, $y_{ 1}$, $y_{13}$, $y_{ 4}$, $y_{15}$, $y_{ 7}$, $y_{14}$\\
$T_{166}$&(2)&
$y_{ 1}$, $y_{ 2}$, $y_{11}$, $y_{15}$, $y_{14}$, $y_{12}$, $y_{ 7}$, $y_{ 6}$, $y_{ 5}$, $y_{ 4}$, $y_{ 8}$, $y_{ 3}$, $y_{13}$, $y_{10}$, $y_{ 9}$\\
$T_{167}$&(2)&
$y_{ 1}$, $y_{ 2}$, $y_{12}$, $y_{ 3}$, $y_{ 5}$, $y_{10}$, $y_{ 9}$, $y_{14}$, $y_{ 7}$, $y_{13}$, $y_{15}$, $y_{ 4}$, $y_{ 8}$, $y_{11}$, $y_{ 6}$\\
$T_{168}$&(1)&
$y_{ 2}$, $y_{ 3}$, $y_{11}$, $y_{ 9}$, $y_{ 8}$, $y_{ 4}$, $y_{ 6}$, $y_{13}$, $y_{10}$, $y_{ 7}$, $y_{15}$, $y_{12}$, $y_{14}$, $y_{ 5}$, $y_{ 1}$\\
$T_{169}$&(1)&
$y_{ 2}$, $y_{ 3}$, $y_{13}$, $y_{ 8}$, $y_{ 5}$, $y_{ 7}$, $y_{ 6}$, $y_{ 9}$, $y_{ 4}$, $y_{11}$, $y_{12}$, $y_{14}$, $y_{ 1}$, $y_{15}$, $y_{10}$\\
$T_{170}$&(1)&
$y_{ 2}$, $y_{ 8}$, $y_{13}$, $y_{ 9}$, $y_{14}$, $y_{12}$, $y_{11}$, $y_{15}$, $y_{10}$, $y_{ 1}$, $y_{ 7}$, $y_{ 3}$, $y_{ 6}$, $y_{ 5}$, $y_{ 4}$\\
$T_{171}$&(2)&
$y_{ 2}$, $y_{11}$, $y_{13}$, $y_{ 5}$, $y_{ 7}$, $y_{ 3}$, $y_{ 1}$, $y_{12}$, $y_{ 4}$, $y_{14}$, $y_{10}$, $y_{ 9}$, $y_{15}$, $y_{ 6}$, $y_{ 8}$\\
$T_{172}$&(2)&
$y_{ 2}$, $y_{ 7}$, $y_{ 9}$, $y_{ 3}$, $y_{ 6}$, $y_{ 8}$, $y_{15}$, $y_{ 1}$, $y_{11}$, $y_{14}$, $y_{ 4}$, $y_{ 5}$, $y_{13}$, $y_{12}$, $y_{10}$\\
$T_{173}$&(1)&
$y_{ 1}$, $y_{ 2}$, $y_{ 4}$, $y_{ 9}$, $y_{15}$, $y_{11}$, $y_{ 7}$, $y_{12}$, $y_{14}$, $y_{10}$, $y_{ 5}$, $y_{ 3}$, $y_{ 6}$, $y_{ 8}$, $y_{13}$\\
$T_{174}$&(2)&
$y_{ 1}$, $y_{ 2}$, $y_{13}$, $y_{ 4}$, $y_{ 7}$, $y_{15}$, $y_{12}$, $y_{ 5}$, $y_{ 6}$, $y_{ 8}$, $y_{11}$, $y_{ 9}$, $y_{10}$, $y_{14}$, $y_{ 3}$\\
$T_{175}$&(1)&
$y_{ 1}$, $y_{ 2}$, $y_{ 4}$, $y_{ 6}$, $y_{10}$, $y_{12}$, $y_{14}$, $y_{ 7}$, $y_{ 3}$, $y_{ 9}$, $y_{11}$, $y_{15}$, $y_{ 8}$, $y_{13}$, $y_{ 5}$\\
$T_{176}$&(2)&
$y_{ 1}$, $y_{ 6}$, $y_{ 7}$, $y_{ 2}$, $y_{ 5}$, $y_{13}$, $y_{ 9}$, $y_{11}$, $y_{ 4}$, $y_{ 3}$, $y_{12}$, $y_{14}$, $y_{15}$, $y_{10}$, $y_{ 8}$\\
$T_{177}$&(2)&
$y_{ 1}$, $y_{ 2}$, $y_{ 6}$, $y_{15}$, $y_{ 9}$, $y_{11}$, $y_{10}$, $y_{12}$, $y_{ 5}$, $y_{ 8}$, $y_{ 3}$, $y_{ 4}$, $y_{13}$, $y_{14}$, $y_{ 7}$\\
$T_{178}$&(1)&
$y_{ 1}$, $y_{ 6}$, $y_{ 8}$, $y_{ 4}$, $y_{ 2}$, $y_{13}$, $y_{ 9}$, $y_{12}$, $y_{11}$, $y_{14}$, $y_{ 7}$, $y_{ 5}$, $y_{15}$, $y_{ 3}$, $y_{10}$\\
$T_{179}$&(1)&
$y_{ 1}$, $y_{ 2}$, $y_{ 3}$, $y_{ 9}$, $y_{15}$, $y_{11}$, $y_{10}$, $y_{ 7}$, $y_{14}$, $y_{12}$, $y_{13}$, $y_{ 8}$, $y_{ 6}$, $y_{ 4}$, $y_{ 5}$\\
$T_{180}$&(1)&
$y_{ 1}$, $y_{ 2}$, $y_{ 6}$, $y_{ 4}$, $y_{10}$, $y_{12}$, $y_{15}$, $y_{14}$, $y_{ 9}$, $y_{ 5}$, $y_{11}$, $y_{ 3}$, $y_{13}$, $y_{ 8}$, $y_{ 7}$\\
$T_{181}$&(1)&
$y_{ 2}$, $y_{10}$, $y_{12}$, $y_{15}$, $y_{ 1}$, $y_{14}$, $y_{ 6}$, $y_{ 8}$, $y_{ 4}$, $y_{ 7}$, $y_{13}$, $y_{ 9}$, $y_{ 3}$, $y_{11}$, $y_{ 5}$\\
$T_{182}$&(1)&
$y_{ 1}$, $y_{ 2}$, $y_{ 6}$, $y_{13}$, $y_{10}$, $y_{11}$, $y_{ 7}$, $y_{ 3}$, $y_{ 9}$, $y_{15}$, $y_{ 8}$, $y_{12}$, $y_{ 4}$, $y_{14}$, $y_{ 5}$\\
$T_{183}$&(2)&
$y_{ 1}$, $y_{ 2}$, $y_{ 3}$, $y_{ 8}$, $y_{ 9}$, $y_{15}$, $y_{ 5}$, $y_{11}$, $y_{10}$, $y_{ 6}$, $y_{ 7}$, $y_{14}$, $y_{ 4}$, $y_{13}$, $y_{12}$\\
$T_{184}$&(1)&
$y_{ 2}$, $y_{ 5}$, $y_{10}$, $y_{15}$, $y_{ 1}$, $y_{14}$, $y_{ 6}$, $y_{12}$, $y_{11}$, $y_{ 7}$, $y_{13}$, $y_{ 4}$, $y_{ 3}$, $y_{ 9}$, $y_{ 8}$\\
$T_{185}$&(1)&
$y_{ 1}$, $y_{ 2}$, $y_{12}$, $y_{ 8}$, $y_{14}$, $y_{15}$, $y_{ 4}$, $y_{ 9}$, $y_{ 5}$, $y_{11}$, $y_{10}$, $y_{13}$, $y_{ 7}$, $y_{ 3}$, $y_{ 6}$\\
$T_{186}$&(1)&
$y_{ 1}$, $y_{ 2}$, $y_{ 4}$, $y_{10}$, $y_{12}$, $y_{ 6}$, $y_{ 7}$, $y_{ 5}$, $y_{ 3}$, $y_{14}$, $y_{11}$, $y_{15}$, $y_{ 8}$, $y_{13}$, $y_{ 9}$\\
$T_{187}$&(1)&
$y_{ 2}$, $y_{14}$, $y_{13}$, $y_{15}$, $y_{ 3}$, $y_{ 4}$, $y_{10}$, $y_{ 8}$, $y_{11}$, $y_{ 1}$, $y_{12}$, $y_{ 5}$, $y_{ 6}$, $y_{ 9}$, $y_{ 7}$\\
$T_{188}$&(2)&
$y_{ 1}$, $y_{ 3}$, $y_{ 4}$, $y_{ 9}$, $y_{14}$, $y_{ 6}$, $y_{15}$, $y_{ 2}$, $y_{11}$, $y_{10}$, $y_{ 5}$, $y_{ 8}$, $y_{12}$, $y_{ 7}$, $y_{13}$\\
$T_{189}$&(1)&
$y_{ 1}$, $y_{ 2}$, $y_{ 7}$, $y_{10}$, $y_{12}$, $y_{ 4}$, $y_{ 6}$, $y_{14}$, $y_{ 3}$, $y_{11}$, $y_{ 8}$, $y_{ 5}$, $y_{ 9}$, $y_{15}$, $y_{13}$\\
$T_{190}$&(2)&
$y_{ 1}$, $y_{ 2}$, $y_{ 4}$, $y_{15}$, $y_{12}$, $y_{ 6}$, $y_{ 7}$, $y_{ 5}$, $y_{ 3}$, $y_{14}$, $y_{11}$, $y_{10}$, $y_{ 8}$, $y_{13}$, $y_{ 9}$\\
$T_{191}$&(2)&
$y_{ 1}$, $y_{ 2}$, $y_{ 3}$, $y_{ 4}$, $y_{14}$, $y_{15}$, $y_{10}$, $y_{ 7}$, $y_{ 6}$, $y_{13}$, $y_{12}$, $y_{11}$, $y_{ 8}$, $y_{ 5}$, $y_{ 9}$\\

      \hline
  \end{tabular}
\end{center}


\newpage

\begin{thebibliography}{00}


\bibitem{Ballmann-Brin} W.~Ballmann and M.~Brin, {\em Polygonal complexes and
    combinatorial group theory}, Geometriae Dedicata 50 (1994),
  165--191.


\bibitem{Bourdon} M.~Bourdon, {\em Immeubles hyperboliques, dimension
    conforme et rigidit\'e de Mostow}, Geom.  Funct. Anal. 7 (1997),
  245--268.

\bibitem{Bourdon2} M.~Bourdon, {\em Sur les immeubles fuchsiens et
    leur type de quasi-isom\'etrie} (French) [Fuchsian buildings and
  their quasi-isometry type] Ergodic Theory Dynam. Systems 20 (2000),
  no. 2, 343--364.

 

\bibitem{CT} I.~Capdeboscq (Korchagina) and A. ~Thomas. {\em Cocompact
    lattices of minimal covolume in rank 2 Kac-Moody groups, Part I:
    Edge-transitive lattices}, arXiv:0907.1350v1 [math.GR]

\bibitem{Carbone-Cartwright-Steger} L.~Carbone, D.~Cartwright and
  T.~Steger, {\em Cocompact lattices in hyperbolic Kac-Moody groups},
  Preprint, (2006)

\bibitem{CarboneCobbs} L.~Carbone and C.~Cobbs, {\em Infinite
    descending chains of cocompact lattices in Kac-Moody groups}, To
  appear in Journal of Algebra and its Applications, (2011)

\bibitem{CarboneGarland} L.~Carbone and H.~Garland,  {\em Existence
    of Lattices in Kac-Moody Groups over Finite Fields},
  Communications in Contemporary Math, Vol 5, No.5, (2003), 813--867.

\bibitem{Cartwright} D.~Cartwright, A.~Mantero, T.~Steger, A.~Zappa,
  {\em Groups acting simply transitively on vertices of a building of
    type $\tilde{A}_2$}, Geometriae Dedicata 47 (1993), 143--166.

\bibitem{Cartwright2} D.~Cartwright, A.~Mantero, T.~Steger and A.~Zappa,
  {\em Groups acting simply transitively on vertices of a building of
    type $\tilde{A}_2$, II The cases $q=2$ and $q=3$}, Geometriae
  Dedicata 47 (1993), 167--223.

\bibitem{Cartwright3} D.~Cartwright and T.~Steger, {\em Enumeration of
    the 50 fake projective planes}, C. R. Math. Acad. Sci. Paris 348
  (2010), no. 1-2.

\bibitem{Edjvet-Howie} M.~Edjvet and J.~Howie, {\em Star graphs,
    projective planes and free subgroups in small cancellation
    groups}, Proc.  London Math. Soc. (3) 57 (1988), no. 2, 301--328.

\bibitem{Gaboriau-Paulin}D.~Gaboriau and F.~Paulin, {\em Sur les
    immeubles hyperboliques}, Geometriae Dedicata 88 (2001), no. 1--3,
  153--197.

\bibitem{Ghys-Harpe} E.~Ghys, P. de la Harpe (eds.), {\em Sur les
    groupes Hyperboliques d'apr\`es Mikhael Gromov}, Birh\"auser, \
  Boston, Basel, \ Berlin, \ 1990.

\bibitem{Haglund} F.~Haglund, {\em Existence, uniqueness and
    homogeneity of certain hyperbolic buildings}, Math. Z.  242
  (2002), no. 1, 97--148.

\bibitem{KV} R.~Kangaslampi and A.~Vdovina, {\em Cocompact actions on
    hyperbolic buildings}, Internat. J. Algebra Comput. 20 (2010),
  no. 4, 591--603.

\bibitem{Kato} F.~Kato and H.~Ochiai, {\em Arithmetic structure of CMSZ
    fake projective planes.}  J. Algebra 305 (2006), no. 2, 116--1185


\bibitem{Jacek} J.~Swiatkowski, {\em Trivalent polygonal complexes of
    nonpositive curvature and Platonic symmetry.}  Geom. Dedicata 70
  (1998), no. 1, 87--110.

\bibitem{Remy} B.~R\'emy,  {\em Groupes de Kac-Moody d\'eploy\'es et
    presque d\'eploy\'es. (French) [Split and almost split Kac-Moody
    groups]} \rm Ast\'erisque No. 277 (2002), viii+348 pp.

\bibitem{RemyRonan} B.~R\'emy and M.~Ronan, {\em Topological groups
    of Kac-Moody type, right-angled twinnings and their lattices. }
  Commentarii Mathematici Helvetici 81 (2006), 191--219.

\bibitem{Tits} J.~Tits and R.~M.~Weiss. {\em Moufang polygons}, Springer
  Monographs in Mathematics, Springer-Verlag, Berlin, 2002.

\bibitem{Vdovina} A.~Vdovina, {\em Combinatorial structure of some
    hyperbolic buildings}, Math. Z. 241 (2002), no. 3, 471--478.

\bibitem{Wise} D.~Wise, {\em The residual finiteness of negatively
    curved polygons of finite groups}, Invent.Math. 149(3): 579 --
  617, 2002.

\bibitem{Xie}X.~Xie, {\em Quasi-isometric rigidity of Fuchsian
    buildings}, Topology 45 (2006), no. 1, 101--169.

\end{thebibliography}
\end{document}